\documentclass[12pt]{amsart}
\usepackage{amssymb,amsmath}
\numberwithin{equation}{section}

\def\dib{\bar\partial}
\def\di{\partial}
\def\simleq{\underset\sim<}
\def\simgeq{\underset\sim>}

\def\T{\text}
\def\1#1{\overline{#1}}
\def\2#1{\widetilde{#1}}
\def\3#1{\widehat{#1}}
\def\4#1{\mathbb{#1}}
\def\5#1{\frak{#1}}
\def\6#1{{\mathcal{#1}}}
\def\C{{\4C}}
\def\R{{\4R}}

\def\Re{{\sf Re}\,}
\def\Im{{\sf Im}\,}

\def\phi{\varphi}

\newtheorem{Thm}{Theorem}[section]
\newtheorem{Cor}[Thm]{Corollary}
\newtheorem{Pro}[Thm]{Proposition}
\newtheorem{Lem}[Thm]{Lemma}
\theoremstyle{definition}\newtheorem{Def}[Thm]{Definition}
\theoremstyle{remark}
\newtheorem{Rem}[Thm]{Remark}
\newtheorem{Exa}[Thm]{Example}

\def\Label#1{\label{#1}}
\def\bl{\begin{Lem}}
\def\el{\end{Lem}}
\def\bp{\begin{Pro}}
\def\ep{\end{Pro}}
\def\bt{\begin{Thm}}
\def\et{\end{Thm}}
\def\bc{\begin{Cor}}
\def\ec{\end{Cor}}
\def\bd{\begin{Def}}
\def\ed{\end{Def}}
\def\br{\begin{Rem}}
\def\er{\end{Rem}}
\def\be{\begin{Exa}}
\def\ee{\end{Exa}}
\def\bpf{\begin{proof}}
\def\epf{\end{proof}}
\def\ben{\begin{enumerate}}
\def\een{\end{enumerate}}

\def\1alpha{[\frac1\alpha]}
\def\T{\text}
\def\R{{\Bbb R}}

\def\C{{\Bbb C}}

\numberwithin{equation}{section}
\def\T{\text}

\newcommand{\no}[1]{\|{#1}\|}
\newcommand{\NO}[1]{{\|{#1}\|}^2}
\newtheorem{theorem}{Theorem  }[section]
\newtheorem{definition}[theorem]{Definition }
\newtheorem{lemma}[theorem]{Lemma  }
\newtheorem{proposition}[theorem]{Proposition  }
\newtheorem{corollary}[theorem]{Corollary }
\newtheorem{example}[theorem]{\it Example }

\begin{document}
\title[Loss of derivatives...]{Loss of derivatives in the infinite type}        
\author[ T.V.~Khanh, S.~Pinton and G.~Zampieri ]
{Tran Vu Khanh, Stefano Pinton and Giuseppe Zampieri}
\address{Dipartimento di Matematica, Universit\`a di Padova, via 
Trieste 63, 35121 Padova, Italy}
\email{khanh@math.unipd.it, pinton@math.unipd.it, 
zampieri@math.unipd.it}
\maketitle

\begin{abstract}
We discuss loss of derivatives for degenerate vector fields obtained from infinite type exponentially non-degenerate hypersurfaces of $\C^2$.

\end{abstract}
\section{Introduction}
\Label{s1}
A system of vector fields  $\{L_j\}$ has subelliptic estimates when it has a gain of  $\delta>0$ derivatives in the sense that $\NO{\Lambda^\delta u}\simleq \sum_j\NO{L_ju}+\NO{u}$, $u\in C^\infty_c$. Here $\Lambda$ is the standard elliptic pseudodifferential operator of order $1$. A system which has finite bracket  type $2m$ is a system whose commutators of order $2m-1$ span the whole tangent space. It is classical that finite type $2m$ implies $\delta$-subelliptic estimates for $\delta=\frac1 {2m}$. If  $\T{Span}\{L,\bar L\}$, in $\C\times\R$, is identified to the tangential bundle $T^{1,0}M\oplus T^{0,1}M$  to a pseudoconvex hypersurface $M\subset \C^2$, then $\{L,\bar L\}$ has finite type $2m$ if and only if the contact of a complex curve $\gamma$ with $M$ is at most $2m$. If the hypersurface is ``rigid", that is, graphed by $\Re w=g(z)$ for a real $C^\infty$ function $g$, then with the notation $g_1=\di_zg,\,\,g_{1\bar1}=\di_z\di_{\bar z}g$ and $t=\Im w$, we have $L=\di_z-ig_1(z)\di_t$ and $[L,\bar L]=g_{1\bar1}\di_t$. It is assumed that $M$ is pseudoconvex, that is, $g_{1\bar1}\geq0$ (this also motivates why the type is  $2m$, even). In terms of $g$, the condition of finite type $2m$ reads
\begin{equation}\Label{1.1}
g_{1\bar1}=0^{2(m-1)}\quad \T{ but }\quad g_{1\bar1}\neq 0^{2m-1}.
\end{equation}
In particular, if $g_{1\bar1}\simgeq|x|^{2(m-1)}$, then we have $\frac1{2m}$-subelliptic estimates.

A system has a superlogarithmic estimate if it has logarithmic gain of  derivative with  an arbitrarily large constant, that is, for any $\delta$ and for suitable $c_\delta$
\begin{equation}
\Label{1.2}
\NO{\log(\Lambda)u}\simleq \delta\sum_j\NO{L_ju}+c_\delta\NO{u},\qquad u\in C^\infty_c.
\end{equation}
A system which satisfies \eqref{1.2} is ``precisely $H^s$-hypoelliptic" for any $s$: $u$ is $H^s$ exactly where the $L_ju$'s are (Kohn \cite{K02}). In particular, the system is $C^\infty$-hypoelliptic. Let $L=\di_z-ig_1(z)\di_t$ for $g$ of infinite type but exponentially non-degenerate in the sense that
\begin{equation}
\Label{1.3}
|z|^\alpha|\log g_{1\bar1}|\searrow 0\T{ as $|z|\searrow 0$ for $\alpha\leq1$}.
\end{equation}
Under this assumption,  $\{L,\bar L\}$ enjoys superlogarithmic estimates (cf. e.g. \cite{KZ10}). If we consider the perturbed system $\{\bar L,\bar z^kL\}$ (any fixed $k\geq 1$), the system has no more superlogarithmic estimates, in general; if $k>1$, a logarithmic loss occurs (Proposition~\ref{p1.1} below).
However, notice that ${\mathcal Lie}\{\bar L,\bar z^kL\}$, the span of commutators of order $\leq k-1$, has superlogarithmic estimates (since it gains $L$). We are able to prove here that  $\{\bar L,\bar z^kL\}$ has, in the terminology of Kohn \cite{K05}, loss of $\frac12$ derivative and thus, in particular, is $C^\infty$-, but not exactly $H^s$-, hypoelliptic. Let $\zeta_0$ and $\zeta_1$ be cut-off functions in a neighborhood of $0$ with $\zeta_0 \prec\zeta_1$ in the sense that $\zeta_1|_{\T{supp}\zeta_0}\equiv1$.
\bt
\Label{t1.1}
Let $L=\di_z-ig_1(z)\di_t$ and assume that $0$ be a point of infinite type, that is, $g_{1\bar1}=0^\infty$ but not exponentially degenerate, that is, \eqref{1.3} be fulfilled. Then the system $\{\bar L,\bar z^kL\}$ (any $k$)  has loss of $\frac12$ derivatives, that is,
\begin{equation}
\Label{1.4}
\NO{\zeta_0u}_s\simleq \NO{\zeta_1\bar Lu}_{s+\frac12}+\NO{\zeta_1\bar z^kLu}_{s+\frac12}+\NO{\bar z^ku}_{\frac12}.
\end{equation}
\et
The proof of this, and the two theorems below, follows in Section~\ref{s2}. 
\br 
\Label{r1.1}
\eqref{1.4} implies local hypoellipticity. Reason is that the loss of derivative takes place only in $t$ (whereas there is elliptic gain in $z$), combined with the fact that $L$ and $\bar L$ have coefficients which are constant in $t$. Thus, if we make a partial regularization $u_\nu\to u,\,\,u_\nu\in C^\infty $ with respect to $t$, use the relation $\bar Lu_\nu=(\bar Lu)_\nu$ (and the same for $L$), and apply \eqref{1.4} to the $u_\nu$'s, we get the proof of the claim.
\er
For $k=1$ we have  estimate for local regularity without  loss
\bt
\Label{t1.2}
In the situation above, assume in addition 
\begin{equation}
\Label{0.0}
|g_1|\simleq g_{1\bar1}^{\frac12};
\end{equation}
 then
\begin{equation}
\Label{0.1}
\NO{\zeta_0u}_s\simleq \NO{\zeta_1\bar Lu}_s+\NO{\zeta_1\bar z Lu}_{s}+\NO{\bar zu}_0.
\end{equation}
\et
When $k>1$, loss may occur
\bp
\Label{p1.1}
Assume that $g=e^{-\frac1{|z|^\alpha}}$. If
\begin{equation}
\Label{0.2}
\NO{\zeta_0u}_s\simleq \NO{(\log\Lambda)^r\zeta_1\bar Lu}_s+\NO{(\log\Lambda)^r\zeta_1\bar z^k Lu}_s+\NO{\bar z^ku}_{\frac12},
\end{equation}
then we must have $r\simgeq \frac{k-(\alpha+1)}\alpha$.
\ep

This seems to be the first time that degenerate vector fields $\{\bar L,\bar z^kL\}$ obtained from $L=\di_z-ig_1(z)\di_t$  of infinite type, that satisfying $g_{1\bar1}=0^\infty$, is considered. However, some additional hypothesis such as  \eqref{1.3}, must be required. This guarantees superlogarithmic estimates (\cite{KZ10}), and in turn, hypoellipticity according to Kohn \cite{K02}.
Loss of derivatives  for $L=\di_z-i\bar z\di_t$ was discovered by Kohn in \cite{K05}. In this case, $L$ is the $(1,0)$ vector field tangential to the strictly pseudoconvex hypersurface $\Re w=|z|^2$ and the loss amounts in $\frac{k-1}2$. The problem was further discussed by Bove, Derridj, Kohn and Tartakoff in \cite{BDKT06} essentially for the vector field $L=\di_z-i\bar z|z|^{2(m-1)}\di_t$ tangential to the hypersurface $\Re w=|z|^{2m}$ and the corresponding loss is $\frac{k-1}{2m}$. In both cases the result extends to the sum of squares $L\bar L+\bar L|z|^{2k}L$ and the loss doubles to $\frac{k-1}m$. For  vector fields $L=\di_z-ig_1(z)\di_t$ tangential to general pseudoconvex  hypersurfaces of finite type (with $g_{1\bar1}$ vanishing at order $2(m-1)$), loss of $\frac{k-1}{2m}$ derivatives has been proved by the authors in \cite{KPZ10}. Under some additional conditions, the result also extends to sums of squares (with doubled loss).

\section{Technical preliminaries and Proof}
\Label{s2}
Our ambient is $\C\times\R$ identified to $\R^3$ with coordinates $(z,\bar z,t)$ or $(\Re z,\Im z, t)$.
We denote by $\xi=(\xi_z,\xi_{\bar z},\xi_t)$ the variables dual to $(z,\bar z,t)$, by $\Lambda^s_\xi$ the standard symbol $(1+|\xi|^2)^{\frac s2}$, and by 
$\Lambda^s$ the pseudodifferential operator with symbol $\Lambda_\xi^s$; this is defined by $\Lambda^s(u)=\mathcal F^{-1}(\Lambda_\xi
^s\mathcal F(u))$ where $\mathcal F$ is the Fourier transform. 
We also consider the partial symbol $\Lambda^s_{\xi_t}$ and the associate pseudodifferential operator $\Lambda^s_t$. 
We denote by $\no{u}_s:=\no{\Lambda^su}_0$ (resp. $\no{u}_{\R, \,s}:=\no{\Lambda_t^su}_0$) the full (resp. totally real) $s$-Sobolev norm.
We use the notation $\simgeq$ and $\simleq$ to denote inequalities up to multiplicative constants; we denote by $\sim$ the combination of $\simgeq$ and $\simleq$.
In $\R^{3}_\xi$, we consider a conical partition of the unity $1=\psi^++\psi^++\psi^0$ where $\psi^\pm$ have support in a neighborhood of the axes $\pm\xi_t$ and $\psi^0$ in a neighborhood of the plane $\xi_t=0$, and introduce a decomposition of the identity $\T{id}=\Psi^++\Psi^-+\Psi^0$ by means of $\Psi^{\overset\pm0}$, the pseudodifferential operators with symbols $\psi^{\overset\pm0}$; we accordingly write $u=u^++u^-+u^0$.
Since $|\xi_z|+|\xi_{\bar z}|\simleq \xi_t$ over $\T{supp}\,\psi^+$, then  $\no{u^+}_{\R, \,s}\sim\no{u^+}_s$.

We carry on the discussion by describing the properties of commutation of the vector fields  $L$ and $\bar L$ for $L=\di_z-ig_1(z)\di_t$.
The crucial equality is 
\begin{equation}
\Label{2.1}
\no{L  u}^2=([L ,\bar L ]u,u)+\no{\bar L  u}^2,\quad u\in C^\infty_c,
\end{equation}
which is readily verified by integration by parts. 
Note here that errors coming from derivatives of coefficients do not occur since $g_1$ does not depend on $t$. 
Since $\sigma(\di_t)$, the symbol of $\di_t$, is dominated by $\sigma(L)$ and $\sigma(\bar L)$ in the ``elliptic region" (the support of $\psi^0$) and since $L$ can be controlled by $\bar L$ with an additional $\epsilon \di_t$ (because of \eqref{2.1}), then $\NO{u^0}_1\simleq \NO{\bar Lu^0}_0+\NO{u}_0$. As for $u^-$, recall that $[L,\bar L]=g_{1\bar1}\di_t$ and hence  $g_{1\bar1}\sigma(\di_t)\leq0$ over $\T{supp}\psi^-$. Thus \eqref{2.1} yields $\NO{Lu}\simleq \NO{\bar Lu}$. It follows that, if $L$ and $\bar L$ have superlogarithmic estimates as in our application, then
$$
\NO{\log(\Lambda)u^-}\leq \delta\NO{\bar Lu^-}+c_\delta\NO{u}.
$$
In conclusion, only estimating $u^+$ is relevant. 
We note here that, over $\T{supp}\,\Psi^+$, we have $g_{1\bar1}\xi_t\geq0$; thus
\begin{equation}
\Label{0.3}
\begin{split}
\NO{g_{1\bar1}u^+}_{\frac12}&=|([L,\bar L]u^+,u^+)|
\\
&\leq \NO{Lu^+}+\NO{\bar Lu^+}.
\end{split}
\end{equation}
Following Kohn \cite{K02}, we introduce  a microlocal modification of $\Lambda^s$, denoted by $R^s$; 
this is the pseudodifferential operator with symbol $R^s_\xi:=(1+|\xi|^2)^{\frac{s\sigma(x)}{2}},\,\,\sigma\in C^\infty_c$; often, what is used is in fact
the partial operator in $t$, $R^s_t$ with symbol $R^s_{\xi_t}$. The relevant property of $R^s$ is
$$
\NO{\Lambda^s\zeta_0u}\simleq \NO{R^s\zeta_0u}+\NO{\zeta_0u}\qquad \T{if $\zeta_0\prec\sigma$}.
$$
Thus, $R^s$ is equivalent to $\Lambda^s$ over functions supported in the region where $\sigma\equiv1$. In addition, $\zeta R^s$ better behaves with 
respect to commutation with $L$; in fact, Jacobi equality yields
\begin{equation}
\Label{log}
[\zeta R^s,L]\sim\dot \zeta R^s+\zeta \log(\Lambda)R^s.
\end{equation}
Thus, on one hand we have the disadvantage of the additional $\log(\Lambda)$ in the second term, but we gain much in the cut-off because
\begin{equation}
\Label{2.2}
\T{ $\dot\zeta R^s$ is of order $0$ if supp$\,\dot\zeta\,\,\cap\,\,$supp$\,\sigma=\emptyset$}.
\end{equation}
Property \eqref{2.2} is  crucial in localizing regularity in presence of superlogarithmic estimates.
\vskip0.3cm
\noindent
{\it Proof of Theorem~\ref{t1.1}.}\hskip0.3cm
As it has already been noticed, it suffices to prove \eqref{1.4} only for $u^+$ and for $\no{\cdot}_{\R,\,\,s}$; thus we write for simplicity $u$ and $\no{\cdot}_s$ but mean $u^+$ and $\no{\cdot}_{\R,\,\,s}$.  Moreover, we  can use a cut-off $\zeta=\zeta(t)$ in $t$ only. In fact, for a cut-off $\zeta=\zeta(z)$ we have $[L ,\zeta(z)]=\dot\zeta$ and $\dot\zeta\equiv0$ at $z=0$. On the other hand, $z^{k}L \sim L $ outside $z=0$ which yields   gain of derivatives, instead of loss. In an estimate we call ``good" a term in the right side (upper bound) and ``absorbable"  a term  which comes as a fraction (small constant or sc) of a formerly encountered term. We take cut-off functions in a neighborhood of $0$: $\zeta_0\prec\sigma\prec\zeta_1\prec\zeta'$; we have for $u\in C^\infty$
\begin{equation}
\Label{2.3}
\begin{split}
\NO{\zeta_0u}_s&=\NO{\zeta_0\zeta_1u}_s
\\
&\leq \NO{R^s\zeta_1u}_0+\NO{[R^s,\zeta_0]\zeta_1u}_0+c\NO{u}_0
\\
&\simleq\NO{R^s\zeta_1u}_0+\NO{u}_0
\\
&\simleq \NO{\zeta'R^s\zeta_1u}_0+\NO{u}_0,
\end{split}
\end{equation}
where the  inequality in the third line follows from interpolation in Sobolev spaces and the last  from $\T{supp}(1-\zeta')\cap\T{supp}\sigma=\emptyset$. We have
\begin{equation}
\Label{2.4}
\begin{split}
\NO{\zeta_0u}_s&\underset{\T{by \eqref{2.3}}}\simleq \underset{(a)}{\underbrace{\NO{R^s\zeta_1u}}}+\NO{u}
\\
&\underset{\T{trivial}}\simleq \underset{(b)}{\underbrace{\NO{\log(\Lambda)R^s\zeta_1u}}}+\NO{u}
\\
&\simleq\NO{\log(\Lambda)(\zeta'R^s\zeta_1)u}+\NO{u}
\\
&\leq \delta\left(\NO{L(\zeta'R^s\zeta_1)u}+\NO{\bar L(\zeta'R^s\zeta_1)u}\right)+c_\delta\NO{u}.
\end{split}
\end{equation}
Here, the inequality in the third line is analogous to the last in \eqref{2.3} in addition to the fact that $[\zeta',\log(\Lambda)]R^s=0(\Lambda^{-1})$; the inequality in the fourth line follows from superlogarithmic estimate.  
Using integration by parts, we estimate the first term in the last line
\begin{equation}
\Label{2.5}
\begin{split}
\NO{L(\zeta'R^s\zeta_1)u}&\simleq \NO{\bar L(\zeta'R^s\zeta_1)u}+\left|([L,\bar L](\zeta'R^s\zeta_1)u,(\zeta'R^s\zeta_1)u)\right|
\\
&\simleq\NO{\bar L(\zeta'R^s\zeta_1)u}+lc \NO{[L,\bar L](\zeta'R^s\zeta_1)u}+\underset{\T{absorbed by (a)}}{\underbrace{sc \NO{R^s\zeta_1u}}}.
\end{split}
\end{equation}
Observe that
\begin{equation}
\Label{2.6}
\begin{split}
\sigma([L,\bar L])=g_{1\bar1}\Lambda^1_\xi&\simleq (g_{1\bar1}^{\frac12}\Lambda^{\frac12}_\xi)|z|^k\Lambda^{\frac12}_\xi
\\
&=\sigma([L,\bar L]^{\frac12})|z|^k\Lambda^{\frac12}_\xi.
\end{split}
\end{equation}
It follows
\begin{equation}
\Label{2.7}
\begin{split}
\NO{[L,\bar L](\zeta'R^s\zeta_1)u}&\leq \NO{[L,\bar L]^{\frac12}\Lambda^{\frac12}(\zeta'R^s\zeta_1)z^ku}
\\
&\leq \NO{Lz^k(\zeta'R^s\zeta_1)u}_{\frac12}+\NO{z^k\bar L(\zeta'R^s\zeta_1)u}_{\frac12}.
\end{split}
\end{equation}
We  wish to first discard the second term in the second line of \eqref{2.7}. For this, we recall Jacobi identity and get
\begin{equation}
\Label{2.8}
\begin{split}
[\bar L,\zeta'R^s\zeta_1]&=[\bar L,\zeta']R^s\zeta_1+\zeta'[\bar L,R^s]\zeta_1+\zeta'R^s[\bar L,\zeta_1]
\\
&\sim \underset{\T{$0$-order  by \eqref{2.2}}}{\underbrace{\dot\zeta'R^s\zeta_1}}+\underset{\T{by \eqref{log}}}{\underbrace{\zeta'\log(\Lambda)R^s\zeta_1}}+\underset{\T{$0$-order  by \eqref{2.2}}}{\underbrace{\zeta'R^s\dot \zeta_1}}.
\end{split}
\end{equation}
Thus we can commutate $z^k\bar L$ with $\zeta'R^s\zeta_1$ in \eqref{2.7} up to an error as described in \eqref{2.8} which yiels
\begin{equation*}
\NO{z^k\bar L(\zeta'R^s\zeta_1)u}_{\frac12}\simleq \NO{(\zeta'R^s\zeta_1)z^k\bar Lu}_{\frac12}+\NO{(\zeta'\log(\Lambda)R^s\zeta_1)z^ku}_{\frac12}+\NO{z^ku}_0.
\end{equation*}
On the other hand, since
$[\zeta',\log(\Lambda)]R^s=0(\Lambda^{-1})$, then
\begin{equation*}
\begin{split}
\NO{(\zeta'\log(\Lambda)R^s\zeta_1)z^ku}_{\frac12}&\simleq \NO{(\log(\Lambda)(\zeta'R^s\zeta_1)z^ku}_{\frac12}+\NO{\zeta_1z^ku}_{-\frac12}
\\
&\underset{\T{suplog estimate}}\simleq\underset{\T{absorbed by $2^{\T{nd}}$ line of \eqref{2.7}}}{\underbrace{\delta\left(\NO{L(\zeta'R^s\zeta_1)z^ku}_{\frac12}+\NO{\bar L(\zeta'R^s\zeta_1)z^ku}_{\frac12}\right)}}+\NO{\zeta_1z^ku}_{-\frac12},
\end{split}
\end{equation*}
where we are using the equality $[\Lambda^{\frac12}_t,L]=0$ as well as $[\Lambda^{\frac12},\log(\Lambda)]=0$. In the same way, using again \eqref{2.8}, we commutate $\bar L$ with $(\zeta'R^s\zeta_1)$ in \eqref{2.4} and \eqref{2.5}. What is left, is to estimate the first term in the last line of \eqref{2.7}. First, 
from Jacobi identity we get
$$
[Lz^k,\zeta'R^s\zeta_1]\sim \T{(0-order)}+z^k\zeta'\log(\Lambda)R^s\zeta_1+\T{(0-order)},
$$
so that we are eventually reduced to estimate $\NO{(\zeta'R^s\zeta_1)Lz^ku}$. This is the most difficult operation. We have (by the trivial identity $[L,z^k]=z^{k-1}$)
$$
\NO{(\zeta'R^s\zeta_1)Lz^ku}_{\frac12}=\underset{\T{good}}{\underbrace{\NO{(\zeta'R^s\zeta_1)z^kLu}_{\frac12}}}+\NO{(\zeta'R^s\zeta_1)z^{k-1}u}_{\frac12}.
$$
Next, 
\begin{equation*}
\begin{split}
\underset{(c)}{\underbrace{\NO{(\zeta'R^s\zeta_1)z^{k-1}u}_{\frac12}}}&=(\underset*{\underbrace{(\zeta'R^s\zeta_1)z^{k-1}u}},(\zeta'R^s\zeta_1)[L,z^k]u)_{\frac12}
\\
&=-(*,(\zeta'R^s\zeta_1)z^kLu)_{\frac12}+(*,(\zeta'R^s\zeta_1)Lz^ku)_{\frac12}.
\end{split}
\end{equation*}
Now,
\begin{equation*}
\begin{cases}
\left|(*,(\zeta'R^s\zeta_1)z^kLu)_{\frac12}\right|\leq sc \NO{*}_{\frac12}+\underset{\T{good}}{\underbrace{\NO{(\zeta'R^s\zeta_1)z^kLu}_{\frac12}}}
\\
\begin{split}
\left|(*,(\zeta'R^s\zeta_1)Lz^ku)_{\frac12}\right|&\leq\left|((\zeta'R^s\zeta_1)\bar Lz^{k-1}u,(\zeta'R^s\zeta_1)z^ku)_{\frac12}\right|
\\&+2\Big|(\underset{\T{absorbed by (c)}}{\underbrace *},\underset{(d)}{\underbrace{[L,(\zeta'R^s\zeta_1)]z^ku}})_{\frac12}\Big|.
\end{split}
\end{cases}
\end{equation*}
We estimate (d). We notice that
\begin{equation}
\Label{2.9}
[L,(\zeta'R^s\zeta_1)]\sim \zeta'\log(\Lambda)R^s\zeta_1+\T{(0-order)}.
\end{equation}
We also remark that
\begin{equation}
\Label{2.10}
\begin{cases}
[\Lambda^{\frac12}\zeta',\log(\Lambda)]R^s=0(\Lambda^{-{\frac12}})\quad(i)
\\
[\zeta',\Lambda^{\frac12}]R^s\sim 0(\Lambda^{-{\frac12}})\quad(ii)
\\
[L,\Lambda^{\frac12}]=0\quad(iii).
\end{cases}
\end{equation}
Hence
\begin{equation}
\Label{2.11}
\begin{split}
\NO{(d)}_{\frac12}&\underset{\T{by \eqref{2.9}}}\simleq \NO{(\zeta'\log(\Lambda)R^s\zeta_1)z^ku}_{\frac12}+\NO{z^ku}_{\frac12}
\\
&\underset{\T{ by \eqref{2.10} (i) and (ii)}}\leq \NO{(\log(\Lambda)\zeta'\Lambda^{\frac12}R^s\zeta_1)z^ku}_0+\NO{z^ku}_{\frac12}+\NO{\zeta_1z^ku}_{-{\frac12}}
\\
&\underset{\T{by suplog estimates}}\leq\delta\left(\NO{L(\zeta'\Lambda^{\frac12}R^s\zeta_1z^ku}+\NO{\bar L(\zeta'\Lambda^{\frac12}R^s\zeta_1)z^ku}\right)+c_\delta\NO{z^ku}_{\frac12}.
\end{split}
\end{equation}
Now, the term with $\delta$ is absorbed by the last term in \eqref{2.7} (after we transform $\Lambda^{\frac12}$ into $\no{\cdot}_{\frac12}$ to fit into \eqref{2.7} and use the fact that $[L\zeta',\Lambda^{\frac12}]\sim\Lambda^{\frac12}$). This concludes the proof of \eqref{1.4}.

\hskip15cm$\Box$
\vskip0.3cm
\noindent
{\it Proof of Theorem~\ref{t1.2}.} \hskip0.2cm
As above, we stay in the positive microlocal cone, the support of $\psi^+$,  and consider only derivatives and cut-off with respect to $t$. From the trivial identity  $[L,z]=1$, 
and from $[L,\zeta_0]=\dot\zeta_0g_1$, we get
\begin{equation*}
\begin{split}
\NO{\zeta_0u}&=([L,z]\zeta_0u,\zeta_0u)
\\
&\simleq \NO{\bar z\zeta_0\bar Lu}_s+\NO{ z\zeta_0Lu}_s+\NO{ zg_1\zeta_1u}_s+sc\NO{\zeta_0u}.
\end{split}
\end{equation*}
Now, the last term is absorbed. 
As for the term before
\begin{equation*}
\begin{split}
\NO{zg_1\zeta_1u}_s&\underset{\T{by \eqref{0.0}}}\leq \NO{g_{1\bar1}^{\frac12}\Lambda^{\frac12}\zeta_1u}_{s-\frac12}
\\
&\underset{\T{by \eqref{0.3}}}\leq \NO{\bar zL\zeta_1u}_{s-\frac12}+\NO{\bar z\zeta_1u}_{s-\frac12}\simleq \NO{\zeta_1\bar zLu}_{s-\frac12}+\NO{\bar z\zeta_2u}_{s-\frac12}\quad\T{for $\zeta_2\succ\zeta_1$}.
\end{split}
\end{equation*}
Now, $\NO{\bar z\zeta_2u}_{s-\frac12}$ is not absorbable by $\NO{zg_1\zeta_1u}_s$, but can be estimated by the $0$-norm using induction over $j$ such that $\frac j2\geq s$.

\hskip15cm $\Box$

\vskip0.3cm
\noindent
{\it Proof of Proposition~\ref{p1.1}.} \hskip0.2cm
As ever, we stay in the positive microlocal cone and take derivatives and cut-off only in $t$. We prove the result for $s$ replaced by $0$ and $\frac12$ replaced by $-\epsilon$. The conclusion for general $s$ follows from the fact that $\di_t$ commutes with $L$ and $\bar L$. We define
\begin{equation*}
v_\lambda=e^{-\lambda(e^{-\frac1{|z|^\alpha}}-it+(e^{-\frac1{|z|^\alpha}}-it)^2)}\qquad\lambda>>0.
\end{equation*}
We denote by $-\lambda A$ the term at exponent and note that $\Re \lambda A=\lambda(e^{-\frac1{|z|^\alpha}}+t^2)$. For $L=\di_z+ig_1(z)\di_t$, we have $\bar L v_\lambda=0$ (which is the key point) and moreover
\begin{equation*}
|\bar z^k Lv_\lambda|\sim |z|^{k-(\alpha+1)}e^{-\lambda(e^{-\frac1{|z|^\alpha}}+t^2)}e^{-\frac1{|z|^\alpha}}.
\end{equation*}
We set
\begin{equation*}
\lambda(e^{-\frac1{|z|^\alpha}},t)=(\theta_1,\frac1{\sqrt\lambda}\theta_2).
\end{equation*}
Under this change we have, over $\T{supp}\,\zeta_0$ and $\T{supp}\,\zeta_1$ which implies $\theta_1<<\lambda$,
$$
|z|^{k-(\alpha+1)}=\frac1{(\log\lambda-\log\theta_1)^{\frac{k-(\alpha+1)}{\alpha}}}.
$$
Hence
we interchange
\begin{equation*}
|\bar z^kLv_\lambda|\dashrightarrow \frac1{(\log \lambda)^{\frac{k-(\alpha+1)}\alpha}}\left(\frac{\theta_1+\theta_2^2}{\left(1-\frac{\log \theta_1}{\log \lambda}\right)^{\frac{k-(\alpha+1)}\alpha}}\right)e^{-(\theta_1+\theta_2^2)}.
\end{equation*}
 Notice that $\theta_1<<\lambda$  and hence, for suitable positive  $c_1$ and $c_2$, we have $c_1<\frac{\theta_1+\theta_2^2}{\left(1-\frac{\log \theta_1}{\log \lambda}\right)^{\frac{k-(\alpha+1)}\alpha}}<c_2$, uniformly over $\lambda$. We also interchange
\begin{equation*}
v_\lambda\dashrightarrow e^{-(\theta_1+\theta_2^2)}.
\end{equation*}
Taking $L^2$ norms yields
\begin{equation*}
 \NO{\bar z^kLv_\lambda}\sim \frac1{(\log\lambda)^{\frac{k-(\alpha+1)}\alpha}}\NO{v_\lambda}.
\end{equation*}
So, the effect on  $L^2$ norm of the action  of $\bar z^kL$ over $v_\lambda$ is comparable to $\frac1{(\log\lambda)^{\frac{k-(\alpha+1)}\alpha}}$. We describe now the effect of the pseudodifferential operator $\log(\Lambda_t)$. We claim that
\begin{equation}
\Label{2.16}
\NO{\log(\Lambda_t)e^{-\lambda t^2}}\sim \log\lambda\NO{e^{-\lambda t^2}}.
\end{equation}
This is a consequence of
\begin{equation}
\Label{2.17}
\log(\Lambda_t)e^{-\lambda t^2}\sim \log\lambda e^{-\lambda t^2}+\Big(\log(\Lambda_{\tilde t})e^{-{\tilde t}^2}\Big)\Big|_{\tilde t=\sqrt\lambda t},
\end{equation}
that we go to prove now. 
Using the coordinate change $\tilde\theta=\sqrt\lambda\theta,\,\,\tilde \xi=\frac\xi{\sqrt\lambda}$, we get
\begin{equation*}
\begin{split}
\int e^{i t\xi}\log(\Lambda_\xi)&\Big(\int e^{-i\xi\theta}e^{-\lambda\theta^2}d\theta\Big)d\xi
\\&=\int e^{it\sqrt\lambda\tilde\xi}\Big(\log(\frac1\lambda+|\tilde\xi|^2)^{\frac12}+\log(\sqrt\lambda)\Big)\Big(\int e^{i\tilde \xi\tilde \theta-
\tilde\theta^2}d\tilde\theta\Big)d\tilde\xi
\\&=\log(\sqrt\lambda) e^{-\lambda t^2}+\Big(\log(\Lambda^\lambda_{\tilde t})e^{-{\tilde t}^2}\Big)\Big|_{\tilde t=\sqrt\lambda t},
\end{split}
\end{equation*}
where $\log(\Lambda^\lambda_{\tilde t})$ is the operator with symbol $\log(\frac1\lambda+|\tilde\xi|^2)^{\frac12}$. This proves \eqref{2.17} and in turn the claim \eqref{2.16}.

We compare now the effect over $v_\lambda$ of $\bar z^kL$ with that of $\log(\Lambda_t)$. If
$$
\NO{\zeta_0v_\lambda}\simleq \NO{\zeta_1 (\log \Lambda_t)^r\bar z^k L v_\lambda}+\NO{v_\lambda}_{-\epsilon},
$$
then, since the right side is estimated from above by
$$
\Big((\log\lambda)^r(\log\lambda)^{-\frac{k-(\alpha+1)}\alpha}+\lambda^{-\epsilon}\Big)\NO{v_\lambda},
$$
we must have that the logarithmic term is not infinitesimal which forces  $r\geq \frac{k-(\alpha+1)}\alpha$.

\hskip15cm$\Box$

\end{document}